\documentclass[11pt,twoside]{article} 
\usepackage[utf8]{inputenc}
\usepackage{amsmath, amsthm, amssymb}
\usepackage[margin=1in]{geometry}
\usepackage{enumerate}
\usepackage{enumitem}
\usepackage{mathrsfs}
\usepackage{mathtools}
\usepackage{tikz-cd}
\usepackage{xcolor}
\usepackage{graphicx}
\usepackage{array}
\usepackage{multirow}
\usepackage{makecell}
\usepackage{diagbox}

\newtheorem{theorem}{Theorem}[section]
\newtheorem{corollary}[theorem]{Corollary}
\newtheorem{lemma}[theorem]{Lemma}
\newtheorem{proposition}[theorem]{Proposition}

\theoremstyle{definition}
\newtheorem{definition}[theorem]{Definition}
\newtheorem*{remark*}{Remark}
\newtheorem{remark}[theorem]{Remark}
\newtheorem{example}[theorem]{Example}
\newtheorem*{notation*}{Notation}

\let\emptyset\varnothing

\usepackage{amsfonts}
\usepackage{makeidx}
\usepackage{graphicx}
\usepackage{physics}
\usepackage{array}
\usepackage{float}

\usepackage{hyperref}
\usepackage{subcaption}
\DeclareMathOperator{\st}{st}
\DeclareMathOperator{\lk}{lk}
\DeclareMathOperator{\lex}{lex}

\hypersetup{colorlinks,
linkcolor=red,
citecolor=red}

\title{The Nevo--Santos--Wilson spheres are shellable}
\author{Yirong Yang\\
\small Department of Mathematics\\
\small University of Washington\\
\small Seattle, WA 98195-4350, USA\\
\small \texttt{yyang1@uw.edu}
}
\date{}

\begin{document}

\maketitle

\begin{abstract}
    Nevo, Santos, and Wilson constructed $2^{\Omega(N^d)}$ combinatorially distinct simplicial $(2d-1)$-spheres with $N$ vertices. We prove that all spheres produced by one of their methods are shellable. Combining this with prior results of Kalai, Lee, and Benedetti and Ziegler, we conclude that for all $D \ge 3$, there are $2^{\Theta(N^{\lceil D/2 \rceil})}$ shellable simplicial $D$-spheres with $N$ vertices.
\end{abstract}

\section{Introduction}\label{sec:intro}

The goal of this paper is to establish the asymptotics of the number of shellable $D$-spheres with $N$ vertices, as $N$ grows to infinity. To achieve this, we show that the spheres produced in \cite[Construction 3]{Nevo} by Nevo, Santos, and Wilson are shellable. 

It follows from Steinitz's theorem (see \cite[Chapter 4]{Ziegler}) that all simplicial $2$-spheres can be realized as the boundary complexes of $3$-polytopes. However, in higher dimensions, there are many more simplicial spheres than the boundaries of polytopes. Let $s(D, N)$ denote the number of combinatorially distinct $D$-spheres with $N$ vertices. For $D \ge 4$, Kalai \cite{Kalai} proved that $s(D, N) \ge 2^{\Omega(N^{\lfloor D/2 \rfloor})}$. Pfeifle and Ziegler \cite{Pfeifle} then complemented Kalai's result by showing $s(3, N) \ge 2^{\Omega(N^{5/4})}$. Later, Nevo, Santos, and Wilson \cite{Nevo} improved the lower bound of $s(D, N)$ for odd $D \ge 3$ to $2^{\Omega(N^{\lceil D/2 \rceil})}$. In constrast to these bounds, we know from works by Goodman and Pollack \cite{Goodman} as well as Alon \cite{Alon} that there are only $2^{\Theta(N \log N)}$ combinatorially distinct $D$-polytopes with $N$ vertices for $D \ge 4$. See also a recent preprint by Padrol, Philippe and Santos \cite{Santos} for the current best lower bound.

An important and related result by Bruggesser and Mani \cite{Mani} is that the boundary complexes of simplicial polytopes are always shellable. This naturally leads to the study of shellable spheres. How many shellable spheres are there? How does this number compare to the number of polytopes? These questions were partially answered by Lee's proof in \cite{Lee} that Kalai's spheres are all shellable. Let $s_{\rm{shell}}(D, N)$ denote the number of shellable $D$-spheres with $N$ vertices. Lee's result implies that 
\begin{equation}\label{eqn:Kalai}
   s_{\rm{shell}}(D, N) \ge 2^{\Omega(N^{\lfloor D/2 \rfloor})}.  
\end{equation}

What about the Nevo--Santos--Wilson spheres? The main result of this paper is:
\begin{theorem}\label{thm:mainthm}
The Nevo--Santos--Wilson spheres in {\rm \cite[Construction 3]{Nevo}} are all shellable.
\end{theorem}
On the other hand, Benedetti and Ziegler \cite{Benedetti} proved that for $D \ge 2$, the number of combinatorially distinct locally constructible (LC) $D$-spheres with $M$ facets grows not faster than $2^{D^2M}$. In addition, they proved that shellable spheres are LC. Meanwhile, the Upper Bound Theorem for simplicial spheres by Stanley \cite{Stanley} asserts that a simplicial $D$-sphere with $N$ vertices has at most $O(N^{\lceil D/2\rceil})$ facets. We make the following observation by combining Theorem \ref{thm:mainthm} with Benedetti and Ziegler's results, the bounds in (\ref{eqn:Kalai}), and the Upper Bound Theorem for simplicial spheres. 
\begin{corollary}\label{cor:maincorollary}
\[
s_{\rm{shell}}(D, N) = 2^{\Theta(N^{\lceil D/2 \rceil})} \text{ for all } D \ge 3.
\]
\end{corollary}

The structure of this paper is as follows. Several key definitions and facts related to Nevo, Santos, and Wilson's construction are provided in Section \ref{sec:prelim}. Sections \ref{sec:shellableball} and \ref{sec:shellablesphere} contain a detailed proof of the shellability of the Nevo--Santos--Wilson spheres. Detailed computations leading to Corollary \ref{cor:maincorollary} can be found at the end of Section \ref{sec:shellablesphere}.

\section{Preliminaries}\label{sec:prelim}

\subsection{Basic definitions}

We start with several essential definitions and notations in preparation for the rest of the paper. 

A \emph{simplicial complex} $\Delta$ on a finite vertex set $V$ is a collection of subsets of $V$ such that if $\sigma \in \Delta$ and $\tau \subseteq \sigma$, then $\tau \in \Delta$. The elements of $\Delta$ are called \emph{faces}. The \emph{dimension} of each face $\sigma$ is $\dim \sigma = |\sigma| - 1$. Conventionally we call the $0$-dimensional faces \emph{vertices}, and the $1$-faces \emph{edges}. The dimension of $\Delta$ is $\dim \Delta = \max \{\dim \sigma: \sigma\in \Delta\}$. We say $\Delta$ is \emph{pure} if its maximal faces with respect to inclusion all have the same dimension. In that case, these maximal faces are called \emph{facets} and faces of one dimension less are called \emph{ridges}. 

The \emph{simplex} on $V$, denoted $\overline{V}$, is the collection of all subsets of $V$. For a face $\sigma \in \Delta$, $\overline{\sigma}$ is the collection of all subsets of $\sigma$. Starting from the next section, we blur the difference between a face $\sigma \in \Delta$ and the simplex $\overline{\sigma} \subseteq \Delta$ and denote both as $\sigma$ by abuse of notation. 

For two integers $n_1, n_2$ such that $n_1 < n_2$, define $[n_1, n_2]:= \{n_1, \dots, n_2\}$ and $[n_1] := \{1, \dots, n_1\}$ when $n_1 > 0$. A \emph{path} of length $n-1$ is a $1$-dimensional pure simplicial complex on the vertex set $\{v_1, \dots, v_n\}$ whose facets are $\{v_i, v_{i+1}\}$ for $i \in [n-1]$. We denote this path as $P(v_1, \dots, v_n)$.

Given a $D$-dimensional simplicial complex $\Delta$, we can associate with $\Delta$ its \emph{geometric realization} $\lVert \Delta \rVert$ as follows. For each maximal face $\sigma \in \Delta$, build a $(|\sigma|-1)$-dimensional \emph{geometric simplex} with vertices labeled by elements in $\sigma$. Glue the simplices in a way that every two simplices are identified along their common (possibly empty) face. We say that $\Delta$ is a \emph{simplicial $D$-sphere} (and respectively, a \emph{simplicial $D$-ball}) if $\lVert \Delta \rVert$ is homeomorphic to a $D$-sphere ($D$-ball). In particular, whenever $\lVert \Delta \rVert$ is homeomorphic to a $D$-manifold, every ridge of $\Delta$ is in at most two facets of $\Delta$. 

Let $\sigma \in \Delta$. The \emph{star} and \emph{link} of $\sigma$ in $\Delta$ are respectively defined to be the subcomplexes
\[
\st_\Delta \sigma = \{\tau \in \Delta: \tau \cup \sigma \in \Delta\}
\text{ and }
\lk_\Delta \sigma = \{\tau \in \Delta: \tau \cap \sigma = \emptyset, \tau \cup \sigma \in \Delta\}.
\]

If $\Delta$ and $\Gamma$ are simplicial complexes on disjoint vertex sets $V$ and $V'$, then the \emph{join} of $\Delta$ and $\Gamma$ is the simplicial complex $\Delta * \Gamma = \{F \cup G: F \in \Delta, \ G \in \Gamma\}$. When one of the the complexes, say $\Gamma$, has only a single vertex $v$, then we call $\Delta * \overline{\{v\}}$ (or simply, $\Delta * v$) the \emph{cone} over $\Delta$ with apex $v$.

Among the many equivalent definitions of shellability, we take the following one for this paper (see for example \cite[Remark 8.3(ii)]{Ziegler}).
\begin{definition}\label{def:shellability}
    A pure $D$-dimensional simplicial complex $\Delta$ is \emph{shellable} if there exists a total order of facets $F_1, \dots, F_n$ of $\Delta$ such that for each $i \in [2, n]$, for every $j < i$, there exists some $m < i$ with the property that $F_i \cap F_m$ is a $(D-1)$-face containing $F_i \cap F_j$. We call such an order a \emph{shelling} of $\Delta$. 
\end{definition}
Every shelling of $\Delta$ induces a shelling of $\st_\Delta \sigma$ and $\lk_\Delta \sigma$. Moreover, the join of two shellable complexes is shellable. 

Simplicial complexes form a subclass of \emph{polyhedral complexes}. A polyhedral complex $C$ is a collection of polytopes such that 
\begin{itemize}
    \item if $P \in C$ and $Q$ is a face of $P$, then $Q \in C$, and 
    \item if $P, Q \in C$, then $P \cap Q$ is a common face of $P$ and $Q$. 
\end{itemize}
For more information about polytopes and shellability, see \cite{Ziegler}. Polyhedral complexes naturally come with a geometric realization. All definitions from the beginning of this section can be adapted to polyhedral complexes. For instance, a polyhedral complex is a {\em polyhedral $D$-sphere} ({\em polyhedral $D$-ball}, respectively) if it is homeomorphic to a $D$-sphere ($D$-ball). 

Given a polyhedral $D$-ball $C$, define its \emph{boundary complex} $\partial C$ to be the subcomplex of $C$ whose facets are the $(D-1)$-faces of $C$ that are contained in exactly one facet of $C$. A \emph{triangulation} of a polyhedral complex $C$ is a simplicial complex $\Delta$ such that 
\begin{itemize}
    \item the geometric realization of $\Delta$ coincides with $C$, and 
    \item every face of $\lVert\Delta\rVert$ is contained in a polytope of $C$. 
\end{itemize}

\subsection{The Nevo--Santos--Wilson spheres} \label{subsec:nswconstruction}

In this section, we present the key facts about \cite[Construction 3]{Nevo} as well as introduce some new definitions and notation. We closely follow Nevo, Santos, and Wilson's paper \cite{Nevo} and state their results without proof. The reader is encouraged to check their paper for more details. 

\begin{remark*}
A word about notation: the construction below is based on the join of $d$ paths, each of length $n-1$. This join is a $(2d-1)$-dimensional complex. We let $D:=2d-1$ and mention right away that each sphere produced in \cite[Construction 3]{Nevo} is $D$-dimensional and has $N= dn+\lceil d(n-1)/(d+2)\rceil +1$ vertices. 
\end{remark*}

Let $B$ be the join of $d$ paths of length $n-1$. Suppose $d\ge 2$ and $n \ge 3$ for nontrivial results. For each $\ell\in [d]$, we denote the $\ell$-th path by $P(1^\ell, \dots, n^\ell)$. Then each $(2d-1)$-simplex $\sigma$ in the join $B$ is of the form 
\[
\sigma = \{i_1^1, i_1+1^1, \dots, i_d^d, i_d+1^d\}.
\]

\begin{notation*}
For simplicity, we drop the parentheses and commas in the expression of simplices from now on. Thus $\sigma$ becomes $i_1^1\ i_1+1^1 \cdots i_d^d\ i_d+1^d$. On the other hand, $\sigma$ can be uniquely represented by a $d$-tuple of indices $(i_1, \dots, i_d) \in [n-1] \times \cdots \times [n-1]$, or $(i_1^1, \dots, i_d^d)$ when we would like to clarify the coordinates with the superscripts. We use these notations interchangeably. Let $\sum \sigma$ denote $i_1 + \cdots + i_d$, the \emph{index sum} of $\sigma$. For the rest of the paper, we reserve the letter $\sigma$ for the $(2d-1)$-simplices of $B$. 
\end{notation*}

\begin{figure}[H]
    \centering
    \includegraphics[width=\linewidth]{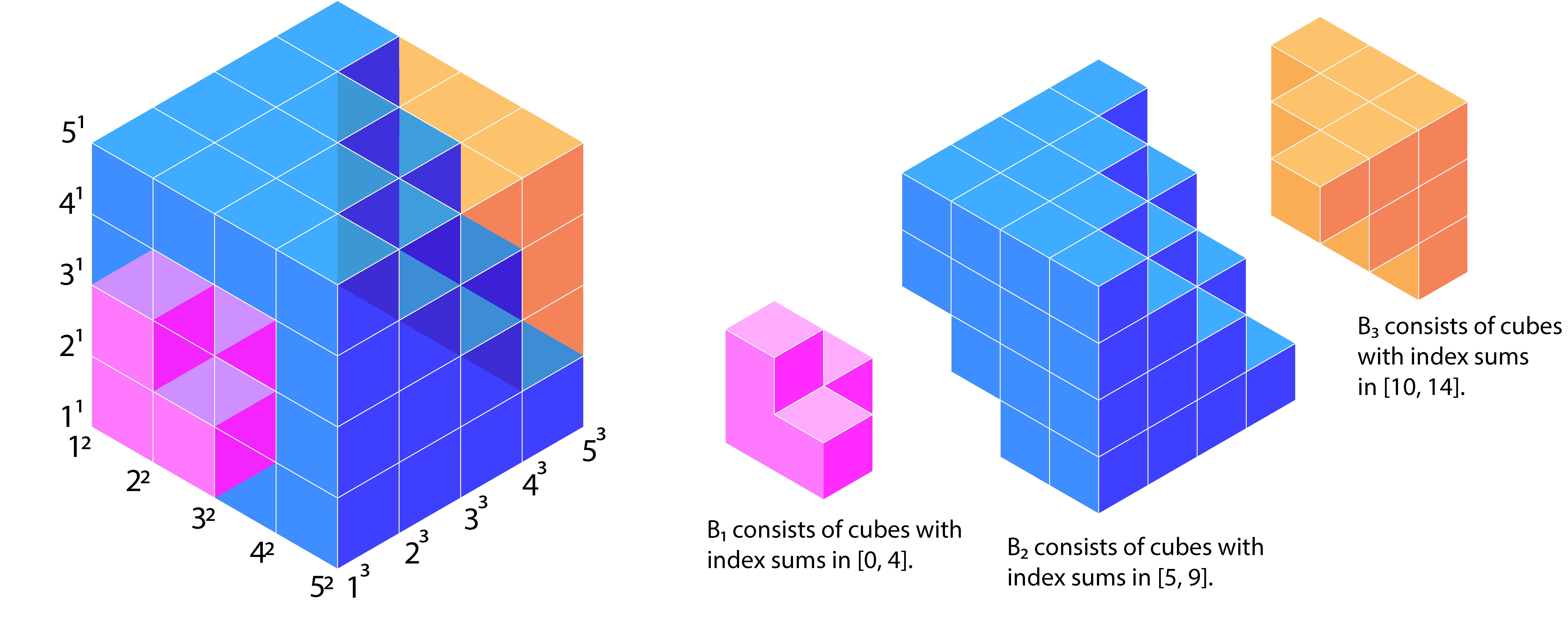}
    \caption{A grid representing $B$ for $d = 3$ and $n = 5$. It is the union of balls $B_k$ for $k \in [3]$.}
    \label{fig:3dgridexample}
\end{figure}

Each $\sigma \in B$ can be visualized as a $d$-cube in the $d$-dimensional grid with side length $n-1$ that represents $B$. The left of Figure \ref{fig:3dgridexample} illustrates the case of $d = 3$ and $n = 5$.
 
The join $B$ is a simplicial $(2d-1)$-ball. For each $k \in [\lceil d(n-1)/(d+2) \rceil\}]$, define $B_k$ to be the union of all $\sigma \in B$ such that $\sum \sigma \in [(k-1) (d+2), k(d+2) -1]$. This is a $(2d-1)$-ball contained in $B$ \cite[Lemma 5.2]{Nevo}. Note that $B = \bigcup_k B_k$. Figure \ref{fig:3dgridexample} depicts the union for the case of $d = 3$ and $n = 5$. 

The idea of \cite[Construction 3]{Nevo} is to replace each $B_k$ with a new ball $\widetilde{B_k}$ with the same boundary. For each $k$, a new vertex $o_k$ is introduced as follows (assuming $B_k$ is not a simplex): 
\begin{itemize}
    \item Consider the collections of $(2d-1)$-simplices:
    \[
    L_k := \left\{\sigma \in B_k: \sum \sigma = (k-1)(d+2)\right\} \text{ and } U_k:= \left\{\sigma \in B_k: \sum \sigma = k(d+2)-1\right\}.
    \]
    We call them respectively the \emph{lower diagonal} and \emph{upper diagonal} of $B_k$. The corresponding cubes for $k = 3$ when $d = 2$ and $n = 8$ are indicated in Figure \ref{fig:2dgridexample}. 
    \begin{figure}[h]
    \centering
    \includegraphics[width=0.6\linewidth]{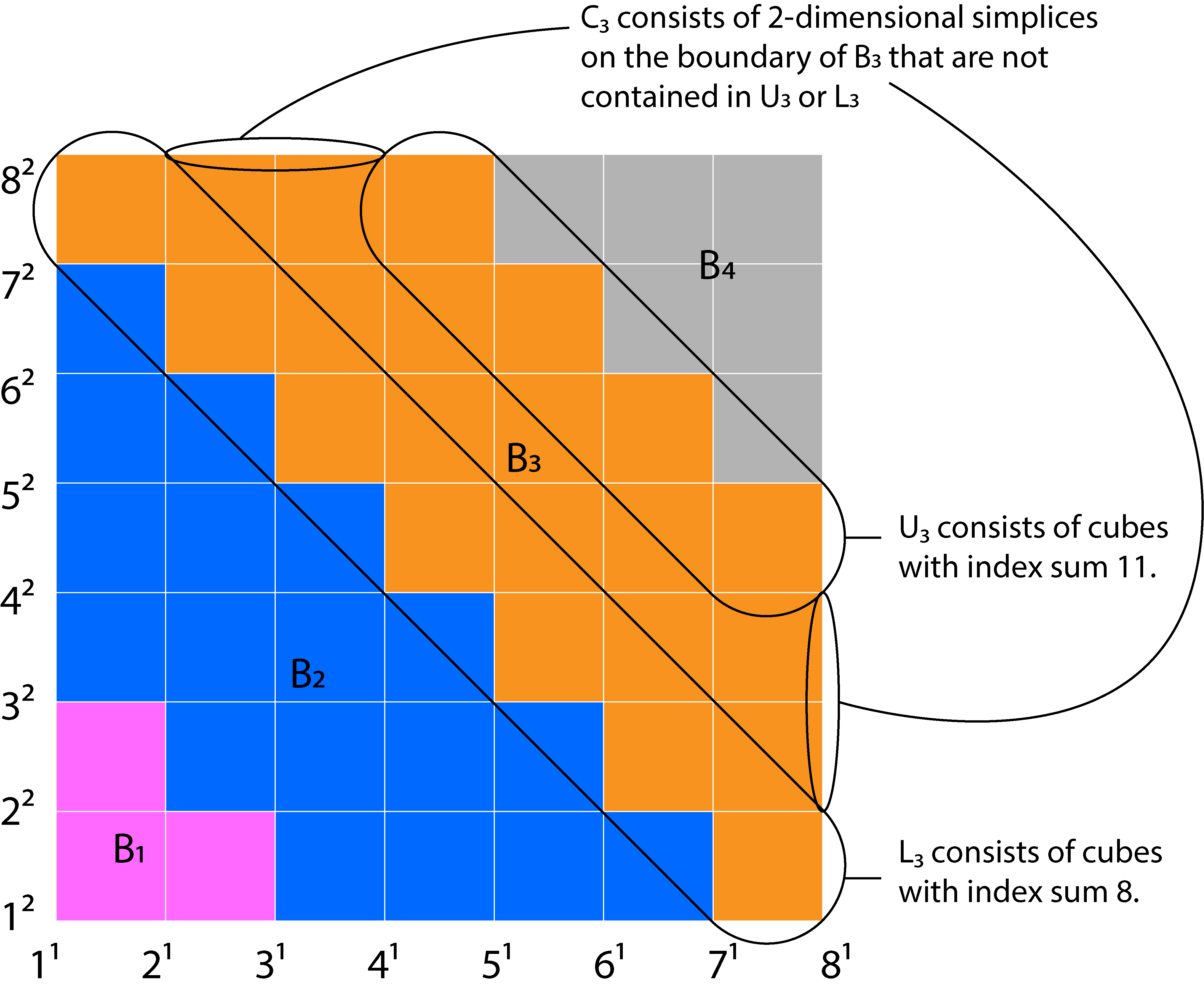}
    \caption{A grid representing $B$ for $d = 2$ and $n = 8$.}
    \label{fig:2dgridexample}
    \end{figure}
    
    For each $\sigma \in L_k \cup U_k$, define $D_\sigma:= \sigma \cap B_k$. Let $C_\sigma$ be a new polyhedral cell whose boundary complex is $D_\sigma \cup (\partial D_\sigma * o_k)$. Let $F_\sigma$ be the only subset of $V(D_\sigma)$ not in $D_\sigma$ but such that all of its proper subsets are in $D_\sigma$. We call $F_\sigma$ the \emph{missing face} of $D_\sigma$. Let $G_\sigma:= \sigma \setminus F_\sigma$. There are two ways to triangulate $C_\sigma$ without introducing new vertices or changing the boundary of $C_\sigma$: 
\[
T_{\sigma, 1} = F_\sigma * \partial(G_\sigma * o_k) \text{ and } T_{\sigma, 2} = \partial F_\sigma * (G_\sigma * o_k).
\]
We use $T_{\sigma, j}$ to mean either one of the two triangulations. 

\item Let $C_k$ be the set of all $(2d-2)$-simplices on the boundary of $B_k$ not contained in any $\sigma \in L_k \cup U_k$. We call $C_k$ the \emph{connecting path} of $B_k$. The corresponding edges for $k = 3$ when $d = 2$ and $n = 8$ are indicated in Figure \ref{fig:2dgridexample}. 
\end{itemize}

Let $\widetilde{B_k} := \bigcup_{\sigma \in L_k \cup U_k} T_{\sigma,j} \cup \bigcup_{\tau \in C_k} \tau * o_k $, and $\widetilde{B} := \bigcup_{k} \widetilde{B_k}$. As the last step, a $(2d-1)$-sphere is obtained by introducing a new vertex $o$ and constructing the complex $\widetilde{B} \cup (\partial \widetilde{B} * o)$. 

As shown in \cite[Corollary 5.5]{Nevo}, this construction yields $2^{\Omega(N^d)}$ combinatorially distinct $(2d-1)$-spheres with $N$ vertices. Indeed, there are two ways to replace $\sigma$ by a triangulation, so there are at least $2^{\sum_k |L_k \cup U_k|} > 2^{2N^d/3d^{d+1}}$ many labeled $N$-vertex triangulations of the $(2d-1)$-sphere. Since $N! = 2^{O(N \log N)}$, dividing by $N!$ does not change the asymptotic order of the bound. 

We end this section by providing a few explicit descriptions about $B$ and the balls $B_k$. These will come in handy for the proofs later. 

\begin{lemma}\label{lem:bboundary}
    The set of facets of $\partial B$ is 
    \[
    \bigcup_{\sigma = (i_1, \dots, i_d) \in B} (\{\sigma \setminus i_\ell+1^\ell: i_\ell = 1\} \cup \{\sigma \setminus i_\ell^\ell: i_\ell = n-1\}).
    \]
\end{lemma}
\begin{proof}
    Any $\tau := \sigma \setminus i_\ell+1^\ell$ is contained in $(i_1^1, \dots, i_\ell - 1^\ell, \dots, i_d^d)$. Therefore, $\tau \in \partial B$ if and only if $i_\ell = 1$. Similarly, any $\tau := \sigma \setminus i_\ell^\ell$ is contained in $(i_1^1, \dots, i_\ell + 1^\ell, \dots, i_d^d)$. Therefore, $\tau \in \partial B$ if and only if $i_\ell = n-1$. 
\end{proof}

\begin{lemma}\label{lem:bkboundary}
The set of facets of $\partial B_k$ is 
\[
\bigcup_{\sigma = (i_1, \dots, i_d) \in L_k} \{\sigma \setminus i_\ell + 1^\ell: \ell \in [d]\} \cup \bigcup_{\sigma = (i_1, \dots, i_d) \in U_k} \{\sigma \setminus i_\ell^\ell: \ell \in [d]\} \cup C_k, 
\]
and
\[
C_k = \bigcup_{\sigma = (i_1, \dots, i_d) \in B_k \setminus (L_k \cup U_k)} (\{\sigma \setminus i_\ell+1^\ell: i_\ell = 1\} \cup \{\sigma \setminus i_\ell^\ell: i_\ell = n-1\}). 
\] 
\end{lemma}
\begin{proof}
The proof is similar to that of Lemma \ref{lem:bkboundary}. Recall that $\sigma \in B_k$ if and only if $\sum \sigma \in [(k-1) (d+2), k(d+2) -1]$. Therefore, $\sigma':= (i_1^1, \dots, i_\ell - 1^\ell, \dots, i_d^d)$ does not belong to $B_k$ if and only if either $i_\ell = 1$ or $\sigma \in L_k$, which forces $\sum \sigma' = \sum \sigma - 1 < (k-1)(d+2)$. Analogously, $\sigma':= (i_1^1, \dots, i_\ell + 1^\ell, \dots, i_d^d)$ does not belong to $B_k$ if and only if either $i_\ell = n-1$ or $\sigma \in U_k$, which forces $\sum \sigma' = \sum \sigma + 1 > k(d+2) -1$. 

This also implies that if $\tau \in C_k$, then since $\tau \subset \sigma$ for some $\sigma$ such that $(k-1) (d+2)< \sum \sigma < k(d+2) -1$, $\tau$ must be of the form $\sigma \setminus i_\ell+1^\ell,\ i_\ell = 1$ or $\sigma \setminus i_\ell^\ell,\ i_\ell = n-1$. 
\end{proof}

\begin{lemma}\label{lem:missingfaces}
    Let $\sigma = (i_1, \dots, i_d) \in L_k \cup U_k$. Then its missing face $F_\sigma$ is 
    \[
    \{i_\ell^\ell: i_\ell = n-1\} \cup \{i_\ell + 1^\ell: \ell \in [d]\}\text{ if } \sigma \in L_k \text{, and } \{i_\ell+1^\ell: i_\ell = 1\} \cup \{i_\ell^\ell: \ell \in [d]\}\text{ if }\sigma \in U_k.
    \]
\end{lemma}
\begin{proof}
    This follows from the definition of $F_\sigma$ and Lemmas \ref{lem:bboundary} and \ref{lem:bkboundary}. 
\end{proof}

\section{The shellability of $\widetilde{B}$}\label{sec:shellableball}

\subsection{Overview} 
To prove that the Nevo--Santos--Wilson spheres from \cite[Construction 3]{Nevo} are shellable, we first prove the shellability of the new ball $\widetilde{B}$ obtained by modifying the join of paths $B$. This is where the major work lies. We start by by outlining the ideas behind and the structure of the proof.

Recall that, to better visualize the join $B$ of $d$ paths of length $n-1$, we identify $B$ with a $d$-dimensional grid with side length $n-1$, where each $d$-cube in the grid represents a $(2d-1)$-simplex in $B$. We can examine the $d$-dimensional grid by ``layers." Each layer is a $(d-1)$-dimensional grid representing a join of $d-1$ paths of length $n-1$: 
\[
B = \bigcup_{i_1 \in [n-1]} i_1^1\ i_1+1^1 * \lk_B i_1^1\ i_1+1^1.
\]
We can further divide each $(d-1)$-dimensional grid $\lk_B i_1^1\ i_1+1^1$ into layers of $(d-2)$-dimensional grids, and so on. Roughly speaking, to modify $B$, we can modify each layer of the grid and stack the new layers together, with some caveats to be addressed soon. This allows us to prove the shellability of $\widetilde{B}$ by induction on $d$. We illustrate the process in the following example before formalizing it. 

\begin{example}\label{ex:shellinginduction}
    Return to the example for $d = 3$ and $n = 5$. On the left of Figure \ref{fig:shellinginduction}, the balls $B_1$, $B_2$, and $B_3$ are colored differently. Each $B_k$ consists of cubes with index sums in $[5(k-1), 5k-1]$. 

    \begin{figure}[h]
        \centering
        \includegraphics[width=1\linewidth]{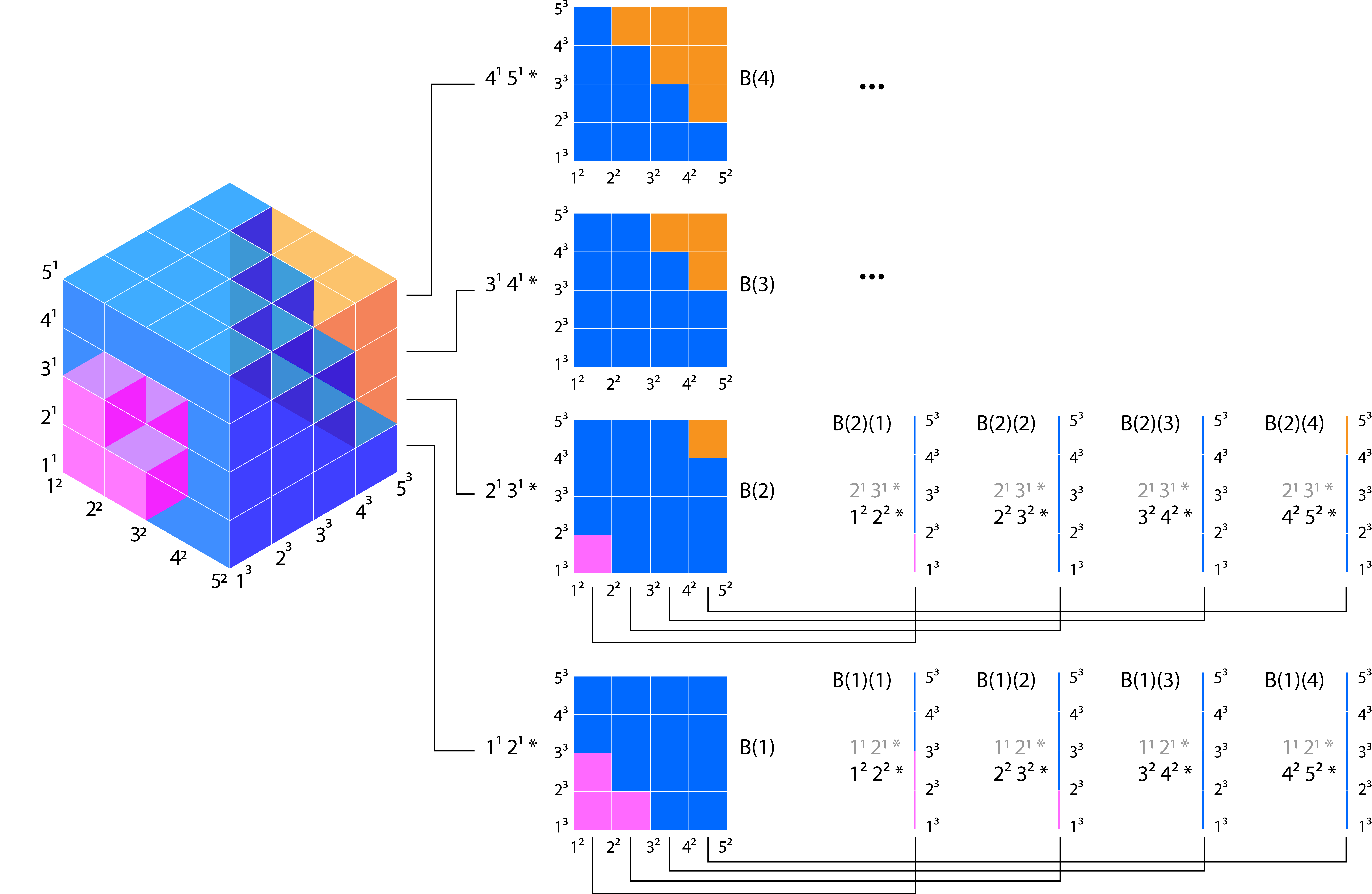}
        \caption{Dividing $B$ into layers recursively.}
        \label{fig:shellinginduction}
    \end{figure}
    
    We divide the $3$-dimensional grid into $4$ layers of $2$-dimensional grids, as shown in the middle column. Let $B(i_1)$ denote the $i_1$-th layer, i.e., $\lk_B i_1^1\ i_1+1^1$. Moreover, let 
    \[
    B(i_1)_k:=\{\tau \in \lk_B i_1^1\ i_1+1^1: i_1^1\ i_1+1^1 \cup \tau \in B_k\}. 
    \]
    Then each $B(i_1)_k$ consists of all $2$-cubes with index sums in $[5(k-1)-i_1, 5k-1-i_1]$. It is colored the same as the corresponding $B_k$ in the $3$-dimensional grid. 
    
    Repeat the process with the $2$-dimensional grids. As shown on the right, we divide the link of $i_1^1\ i_1+1^1$ into $4$ layers of $1$-dimensional grids $B(i_1)(i_2)$. We have reached the base case. Each $B(i_1)(i_2)_k$ consists of all $1$-cubes with index sums in $[5(k-1)-i_1 - i_2, 5k-1-i_1 - i_2]$. 

    Let $\widetilde{B(i_1)}:= \lk_{\widetilde{B}} i_1^1 \ i_1+1^1$. Later, we will show that $\widetilde{B(i_1)}$ can be obtained from $B(i_1)$ following the same rules of modifying $B$. Assuming each $\widetilde{B(i_1)}$ is shellable, then $i_1^1 \ i_1+1^1 * \widetilde{B(i_1)}$ is also shellable. The plan is to concatenate the shellings of $i_1^1 \ i_1+1^1 * \widetilde{B(i_1)}$ together: start with a shelling of $1^1 2^1 * \widetilde{B(1)}$, then attach a shelling of $2^1 3^1 * \widetilde{B(2)}$, and so on. 

    There is one problem: some facets of $\widetilde{B}$ does not belong to any $i_1^1 \ i_1+1^1 * \widetilde{B(i_1)}$. For example, consider $\sigma = (2, 2, 1) = 2^1 3^1 2^2 3^2 1^3 2^3 \in L_2$. Suppose it is replaced by the triangulation $T_{\sigma, 1}$, then one resulting facet is $\sigma \setminus 2^1 \cup o_2$. Another example is $1^1 1^2 2^2 1^3 2^3 o_1$, with $1^1 1^2 2^2 1^3 2^3 \in C_1$. We will characterize these extra facets in Lemma \ref{lem:paddings} and insert them between the shellings of $i_1^1 \ i_1+1^1 * \widetilde{B(i_1)}$.
\end{example}

Let $d \ge 2$, $n \ge 3$, and $B$ be the join of $d$ paths of length $n-1$. We extend the notation in Example \ref{ex:shellinginduction} to the general case. Let $m \in [d]$. Consider the join of $m$ paths 
\[
B(i_1) \cdots (i_{d-m}) := \lk_B (i_1^1\ i_1+1^1  \cdots  i_{d-m}^{d-m} \ i_{d-m}+1^{d-m}). 
\]
When $m = d$, this is just $B$. For $k \in [\lceil d(n-1)/(d+2) \rceil]$, let $B(i_1) \cdots (i_{d-m})_k$ denote the union of all $m$-cubes $(i_{d-m+1}, \dots, i_d)$ with index sums in 
\[
I_m := \left[(k-1) (d+2) - \sum_{\ell = 1}^{d-m} i_\ell,\ k(d+2) -1 - \sum_{\ell = 1}^{d-m} i_\ell\right].
\]
Then $B(i_1) \cdots (i_{d-m}) = \bigcup_k B(i_1) \cdots (i_{d-m})_k$. Just like in $B$, the lower and upper diagonal cubes of $B(i_1) \cdots (i_{d-m})_k$ are respectively those whose index sums attain the minimum and maximum of $I_m$. The connecting path consists of the $(2m-2)$-simplices on the boundary of $B(i_1) \cdots (i_{d-m})_k$ that are not contained in any lower or upper diagonal cubes. Lastly, let
\[
\widetilde{B(i_1) \cdots (i_{d-m})}:= \lk_{\widetilde{B}} (i_1^1\ i_1+1^1 \cdots i_{d-m}^{d-m} \ i_{d-m}+1^{d-m}). 
\]

Based on the inductive idea from Example \ref{ex:shellinginduction}, we would like to prove the statement below, which directly implies the shellability of $\widetilde{B}$.
\begin{proposition}
    Let $m \in [d]$. Then $\widetilde{B(i_1) \cdots (i_{d-m})}$ is shellable for all $i_1, \dots, i_{d-m} \in [n-1]$. 
\end{proposition}

The base case of $m = 1$ is immediately true. Since $\widetilde{B}$ is a ball, the link $\widetilde{B(i_1) \cdots (i_{d-1})}$ is a connected $1$-dimensional simplicial complex, which is automatically shellable. 

\begin{remark}\label{rem:reduce}
For the inductive step, we need to prove the following: let $m \in [2, d]$ and $i_1, \dots, i_{d-m} \in [n-1]$, assuming $\widetilde{B(i_1) \cdots (i_{d-m+1})}$ is shellable for all $i_{d-m+1} \in [n-1]$, then $\widetilde{B(i_1) \cdots (i_{d-m})}$ is shellable. For cleanliness of the write-up, we provide the proof for the shellability of $\widetilde{B}$ assuming $\widetilde{B(i_1)}$ is shellable for all $i_1 \in [n-1]$. The proof for the general case is analogous. The fact that $I_m$ is ``shifted" does not matter in the proof; all we need is that the difference of the endpoints of the interval is greater than the dimension of the grid. To accommodate this, we use the relaxed assumption that the difference between the index sum of the upper diagonal cubes and that of the lower diagonal cubes in $B_k$ is greater than $d$ (instead of being equal to $d+1$). 
\end{remark}

In short, we aim to prove the next proposition in the rest of Section \ref{sec:shellableball}.
\begin{proposition}\label{prop:mainprop}
    Suppose $\widetilde{B(i_1)}$ is shellable for all $i_1 \in [n-1]$. Then $\widetilde{B}$ is shellable. 
\end{proposition}

We claim that the order of facets below is a shelling of $\widetilde{B}$. 

\begin{table}[H]
\setlength{\arrayrulewidth}{0.3mm}
\renewcommand{\arraystretch}{1.6}
\begin{center}
\begin{tabular}{l}
\hline
Start with any shelling order of $1^1 2^1 * \widetilde{B(1)}$.\\
Then list in any order facets of the form $\sigma \setminus 2^1 \cup o_k$ for $\sigma \in L_k \cup U_k \cap \st_B 1^1 2^1$ for all $k$.\\
Then list in any order facets of the form $\sigma \setminus 1^1 \cup o_k$ for $\sigma \in L_k \cup U_k \cap \st_B 1^1 2^1$ for all $k$.\\
\textbf{\emph{Repeat for all $k \in [\lceil d(n-1)/(d+2) \rceil]$ as $k$ increases}}:\\
List the facets of the form $\sigma \setminus 2^1 \cup o_k$ for $\sigma \in (B_k \setminus (L_k \cup U_k)) \cap \st_B 1^1 2^1$\\
in the order $\prec$ defined in Section \ref{subsec:connecting}.\\
\hline
\textbf{\emph{Repeat for all $2 \le i_1 \le n-1$} as $i_1$ increases}:\\
List in any order facets of the form $\sigma \setminus i_1+1^1 \cup o_k$ for $\sigma \in L_k \cup U_k \cap \st_B i_1^1\ i_1+1^1$ for all $k$.\\
Follow this by any shelling order of $i_1^1\ i_1+1^1 * \widetilde{B(i_1)}$.\\
Then list in any order facets of the form $\sigma \setminus i_1^1 \cup o_k$ for $\sigma \in L_k \cup U_k \cap \st_B i_1^1\ i_1+1^1$ for all $k$.\\
\hline
\textbf{\emph{Repeat for all $k \in [\lceil d(n-1)/(d+2) \rceil]$ as $k$ increases}}:\\
List the facets of the form $\sigma \setminus n-1^1 \cup o_k$ for $\sigma \in (B_k \setminus (L_k \cup U_k)) \cap \st_B n-1^1 n^1$\\
in the order $\prec$ defined in Section \ref{subsec:connecting}.\\
\hline
\end{tabular}
\caption{A shelling of $\widetilde{B}$.}
\label{tab:shellingorder}
\end{center}
\end{table}
Here, we do not claim that $\widetilde{B}$ must contain $\sigma \setminus i_1+1^1 \cup o_k$ or $\sigma \setminus i_1^1 \cup o_k$ for every $k$ and every $\sigma \in L_k \cup U_k \cap \st_B i_1^1\ i_1+1^1$. In fact, whether it is a facet depends on the choice of $T_{\sigma, j}$. We merely say that if it is a facet of $\widetilde{B}$, it should be listed in the corresponding slot in this order. 

In Section \ref{subsec:layers}, we analyze the facets of $\widetilde{B}$ that are not contained in $\bigcup_{i_1 \in [n-1]} i_1^1\ i_1+1^1*\widetilde{B(i_1)}$, showing that Table \ref{tab:shellingorder} is indeed the list of all facets of $\widetilde{B}$. In Section \ref{subsec:connecting}, we define a certain order $\prec$ for a specific subset of these facets. We leave the actual proof of the fact that Table \ref{tab:shellingorder} provides a shelling of $\widetilde{B}$ to Section \ref{subsec:minimalnewface}.

\subsection{Facets between layers}\label{subsec:layers}

We first establish that, essentially, the links $B(i_1)$ are modified the same way as $B$ is. 
\begin{lemma}\label{lem:linkoftriangulation}
    Let $i_1 \in [n-1]$. Then $\widetilde{B(i_1)}$ can be obtained from $B(i_1)$ following the same steps of changing $B$ to $\widetilde{B}$ in {\rm \cite[Construction 3]{Nevo}}. 
\end{lemma}
More precisely, if $\sigma \in \st_B i_1^1\ i_1+1^1$ is replaced by the triangulation $T_{\sigma, j}$, then $\sigma \setminus i_1^1\ i_1+1^1 \in B(i_1)$ is replaced by $T_{\sigma \setminus i_1^1\ i_1+1^1, j}$. If a $(2d-2)$-simplex $\tau \in C_k \cap \st_B i_1^1\ i_1+1^1$ is coned with $o_k$, then $\tau \setminus i_1^1\ i_1+1^1$ in the connecting path of $B(i_1)_k$ is coned with $o_k$. 
\begin{proof}
    Let $\sigma \in L_k \cup U_k \cap \st_B i_1^1\ i_1+1^1$. Recall from Section \ref{subsec:nswconstruction} that $D_\sigma = \sigma \cap \partial B_k$, $F_\sigma$ is the missing face of $D_\sigma$, and $G_\sigma = \sigma \setminus F_\sigma$. We extend these notations to $\sigma \setminus i_1^1\ i_1+1^1$. By simple set-theoretic computations, one can check that 
    \begin{align*}
        D_{\sigma \setminus i_1^1\ i_1+1^1} &= (\sigma \setminus i_1^1\ i_1+1^1) \cap \partial B(i_1)_k = \lk_{D_\sigma} i_1^1\ i_1+1^1,\\
        F_{\sigma \setminus i_1^1\ i_1+1^1} &= F_\sigma \setminus i_1^1\ i_1+1^1, \text{ and } G_{\sigma \setminus i_1^1\ i_1+1^1} = G_\sigma \setminus i_1^1\ i_1+1^1.
    \end{align*}
    Let $X = F_\sigma \cap i_1^1\ i_1+1^1$ and $Y = i_1^1\ i_1+1^1 \setminus X$. Then $F_{\sigma \setminus i_1^1\ i_1+1^1} = \lk_{F_\sigma}X$ and $G_{\sigma \setminus i_1^1\ i_1+1^1} = \lk_{G_\sigma}Y$. Observe that the triangulation
    \begin{align*}
        T_{\sigma \setminus i_1^1\ i_1+1^1, 1} &=F_{\sigma \setminus i_1^1\ i_1+1^1} * \partial (G_{\sigma \setminus i_1^1\ i_1+1^1} * o_k)\\
        &= \lk_{F_\sigma}X * \partial (\lk_{G_\sigma}Y * o_k) = \lk_{F_\sigma}X * ( \partial \lk_{G_\sigma}Y * o_k \cup \lk_{G_\sigma}Y *\partial o_k)\\
        &= \lk_{F_\sigma}X * \lk_{\partial G_\sigma}Y * o_k \cup \lk_{F_\sigma}X *\lk_{G_\sigma}Y = \lk_{F_\sigma * \partial G_{\sigma}*o_k} (X \sqcup Y) \cup \lk_{F_\sigma*G_\sigma} (X \sqcup Y)\\
        &= \lk_{T_{\sigma, 1}} i_1^1\ i_1+1^1. 
    \end{align*}
    Similarly, $T_{\sigma \setminus i_1^1\ i_1+1^1, 2} = \lk_{T_{\sigma, 2}} i_1^1\ i_1+1^1$.

    Finally, for $\tau \in C_k \cap \st_B i_1^1\ i_1+1^1$, $\tau \setminus i_1^1\ i_1+1^1 * o_k = \lk_{\tau*o_k} i_1^1\ i_1+1^1$. 
\end{proof}

With this, we can characterize the facets of $\widetilde{B}$ that are ``between" the layers. 

\begin{lemma}\label{lem:paddings}
Let $F$ be a facet of $\widetilde{B} \setminus \bigcup_{i_1 \in [n-1]} i_1^1\ i_1+1^1*\widetilde{B(i_1)}$ that arises from $\sigma = (i_1, \dots, i_d)\in B_k \cap \st_B i_1^1\ i_1+1^1$ for some $k$. 

\begin{enumerate}
    \item If $\sigma \in L_k\cup U_k \cap \st_B i_1^1\ i_1+1^1$, then $F$ takes the following forms depending on $F_\sigma$ and the triangulation $T_{\sigma, j}$:
\setlength{\arrayrulewidth}{0.3mm}
\renewcommand{\arraystretch}{1.6}
\renewcommand\cellset{\renewcommand\arraystretch{1.6}}
\begin{table}[H]
\small
\begin{center}
\begin{tabular}{l|l|l|l|l}
\hline
\diagbox[width=11em]{$T_{\sigma, j}$}{$F_\sigma \cap i_1^1\ i_1+1^1$} & $i_1^1$ & $i_1+1^1$& $\emptyset$ & $i_1^1\ i_1+1^1$\\ \hline
$T_{\sigma, 1}$ & $\sigma \setminus i_1+1^1\cup o_k$ & $\sigma \setminus i_1^1\cup o_k$  & \makecell[l]{$\sigma \setminus i_1+1^1\cup o_k$,\\$\sigma \setminus i_1^1\cup o_k$} & Does not exist  \\ \hline
$T_{\sigma, 2}$ & $\sigma \setminus i_1^1\cup o_k$ & $\sigma \setminus i_1+1^1\cup o_k$  & Does not exist   & \makecell[l]{$\sigma \setminus i_1+1^1\cup o_k$,\\$\sigma \setminus i_1^1\cup o_k$}  \\ \hline
\end{tabular}
\end{center}
\end{table}
\item If $\sigma \in (B_k \setminus (L_k \cup U_k)) \cap \st_B i_1^1\ i_1+1^1$, then $F = \sigma \setminus i_1+1^1 \cup o_k$ when $i_1 = 1$ and $F = \sigma \setminus i_1^1 \cup o_k$ when $i_1 = n-1$.
\end{enumerate}
\end{lemma}
\begin{proof}
    Case 1 is immediate from Lemma \ref{lem:linkoftriangulation}. The facets of $T_{\sigma, j} \setminus \left(i_1^1\ i_1+1^1 * T_{\sigma \setminus i_1^1\ i_1+1^1, j}\right)$ are those with exactly one of $i_1^1$ and $i_1+1^1$. The specific cases follow from the formulas $T_{\sigma, 1} = F_\sigma * \partial(G_\sigma * o_k)$ and $T_{\sigma, 2} = \partial F_\sigma * (G_\sigma * o_k)$. 

    For case 2, by Lemma \ref{lem:bkboundary}, $C_k = \bigcup_{\sigma \in B_k \setminus (L_k \cup U_k)} (\{\sigma \setminus i_\ell+1^\ell: i_\ell = 1\} \cup \{\sigma \setminus i_\ell^\ell: i_\ell = n-1\})$. Therefore, the only facets in $\bigcup_{\tau \in C_k} \tau * o_k$ that are not covered by $\bigcup_{i_1 \in [n-1]} i_1^1\ i_1+1^1* B(i_1)$ are those listed above.
\end{proof}

\subsection{The total order $\prec$}\label{subsec:connecting}

The goal of this section is to define an order $\prec$ for each of 
\begin{align*}
    &\{\sigma \setminus 2^1 \cup o_k: \sigma \in B_k \setminus (L_k \cup U_k) \cap \st_B 1^1 2^1\}\\ =& \left\{1^1 \cup \tau \cup o_k: \tau \in \lk_B 1^1 2^1, \sum \tau \in [(k-1) (d+2),\ k(d+2)-3]\right\}\\
    &\hspace{0.2\textwidth}\text{and}\\
    &\{\sigma \setminus n-1^1 \cup o_k: \sigma \in B_k \setminus (L_k \cup U_k) \cap \st_B n-1^1 n^1\}\\
    =& \left\{n^1 \cup \tau \cup o_k: \tau \in \lk_B n-1^1 n^1, \sum \tau \in [(k-1) (d+2) +1 - n,\ k(d+2) -2 - n]\right\}.
\end{align*}

Let $n_1, n_2$ be positive integers. We consider the \emph{lexicographic order} on the elements of $[n_1]^{n_2}$. For $A, B \in [n_1]^{n_2}$, define $A <_{\lex} B$ if the leftmost nonzero coordinate of $B-A$ is positive. 

\begin{notation*}
    For a set $S$ and a total order $<$ on $S$, let $S_<$ denote the ordered set of elements in $S$ with respect to $<$, and let $S_{<^{-1}}$ denote the ordered set of elements in $S$ with respect to $<$ in reverse. 
\end{notation*}

Consider $T:= \{\tau = (i_2, \dots, i_d): \sum \tau \in [(k-1) (d+2) +1 - a,\ k(d+2) - 2 - a]\}$ for a fixed $k$ and $a \in \{1, n\}$. We define a total order $\prec$ on $T$, which can then be extended to the two sets above. We begin by visualizing this order through an example. 

\begin{example}\label{ex:connectinginduction}
    Consider an example for $d = 4, n = 6, k = 2, a = 1$. The pile of $3$-cubes at the top left corner of Figure \ref{fig:connectinginduction} is $B(1)_2$. The gray cubes are from the lower and upper diagonals of $B(1)_2$. The cubes in $T = \{\tau = (i_2, \dots, i_d): \sum \tau \in [6, 9]\}$ are transparent with black outlines. 

    The order $\prec$ on $T$ needs to be somewhat involved for the following reason. We plan to insert this order (having added $1^1$ and $o_k$ to each $\tau$) after a shelling of $1^1 2^1 * \widetilde{B(1)}$ and some extra facets. Although we are not shelling the gray cubes but their related triangulations, we can still visualize the rough locations of the triangulations by looking at these original cubes. Suppose the gray cubes are already in place as part of a bigger complex, we want to fill in the transparent cubes in a way that each new cube intersects the existing complex nicely. Therefore, although $T$ is in fact shellable, merely inserting a shelling order of $T$ would not suffice. 

    The top right corner shows $T$ only. First, put $\{\tau: \sum \tau = 9\}$ in lexicographic order. The $19$ cubes are in blue and labeled by this order. We use the idea of ``layers" again and construct the rest of the order recursively. 

    \begin{figure}[h]
        \centering
        \includegraphics[width=1\linewidth]{connectinginduction.png}
        \caption{Order $\prec$ on $T = \{\tau \in \lk_B 1^1 2^1: \sum \tau \in [6, 9]\}$.}
        \label{fig:connectinginduction}
    \end{figure}
    
    On the second row, $T$ is divided into $5$ layers, $T \cap \st_{B(1)} i_2^2\ i_2+1^2$ for $i_2 \in [5]$. In each $2$-dimensional layer, the cubes already ordered ($1$ through $19$) are now in gray. Starting from $T \cap \st_{B(1)} 1^2 2^2$, arrange $\{\tau: \sum \tau = 8\}$ in lexicographic order. The corresponding cubes are now in blue and labeled $20$ through $23$.

    Next, divide $T \cap \st_{B(1)} 1^2 2^2$ into $T \cap \st_{B(1)} 1^2 2^2i_3^3\ i_3+1^3$ for $i_3 \in [5]$, shown on the last row. The cubes $1$ through $23$ are now gray. Starting from $T \cap \st_{B(1)} 1^2 2^2 1^3 2^3$, put $\{\tau: \sum \tau = 7\}$ in lexicographic order. There is only one such cube, which we label as $24$. Then we add the last cube in that layer to our order. Repeat the procedure for all other $T \cap \st_{B(1)} 1^2 2^2i_3^3\ i_3+1^3$. After this, return to work on $T \cap \st_{B(1)} 2^2 3^2$: first order the cubes with index sum $8$, then divide it into layers $T \cap \st_{B(1)} 2^2 3^2i_3^3\ i_3+1^3$, and so on. 

    We can view the picture as a tree rooted at $T$. The order is given by traversing the tree from top to bottom and left to right.
\end{example}

We now formalize this order and provide a simple algorithm to compare any two cubes in $T$. We first address the weakened assumption given in Remark \ref{rem:reduce}. That is, the difference between the index sum of the upper diagonal cubes and that of the lower diagonal cubes in $B_k$ is greater than $d$. Therefore, we assume $T = \{\tau: \sum \tau \in [t, t']\}$ for some $t, t'$ such that $t' - t \ge d - 1$. 

For $\tau \in T$, define $I(\tau):= t + d - \sum \tau$. Then for $\tau \in T$ such that $\sum \tau - t \ge d - 1$, $I(\tau) \le 1$. For $\tau \in T$ such that $\sum \tau - t < d - 1$, $I(\tau) \in [2, d]$. In particular, the smaller $I(\tau)$ is, the closer $\tau$ is to the ``upper diagonal" of $T$. 

Let $i_1 \in \{1, n-1\}$. For $m \in [2, d]$ and $i_2, \dots, i_m \in [n-1]$, define $T_m(i_2, \dots, i_m)$ to be an ordered set of cubes in $\{\tau \in T: I(\tau) \in [2, d]\} \cap \st_{B(i_1)}i_2^2\ i_2+1^2 \cdots i_m^m\ i_m+1^m$ as follows: 
\begin{itemize}
    \item If $m = d$, then $T_m(i_2, \dots, i_m)$ either contains only $\tau = i_2^2\ i_2+1^2 \cdots i_m^m\ i_m+1^m$ when $I(\tau) = m = d$, or is empty when $I(\tau) \ne d$. 
    \item If $m < d$, then let 
    \begin{align*}
        T_m(i_2, \dots, i_m) &:= (\{\tau \in T: I(\tau) = m\} \cap \st_{B(i_1)} i_2^2\ i_2+1^2 \cdots i_m^m\ i_m+1^m)_{<_{\lex}}\\
        &\cup T_{m+1}(i_2, \dots, i_m, 1) \cup \cdots \cup T_{m+1}(i_2, \dots, i_m, n-1).
    \end{align*}
\end{itemize}
We concatenate the ordered sets together:
\[
T_\prec: = \{\tau \in T: I(\tau) \le 1\}_{<_{\lex}^{-1}} \cup T_2(1) \cup \cdots \cup T_2(n-1). 
\]

The lemma below provides an algorithm to compare two cubes in $\{\tau \in T: I(\tau) \in [2,d]\}$.
\begin{lemma}\label{lem:connectingalgorithm}
    The order on $T(1, 2) \cup \cdots \cup T(n-1, 2)$ defined above is equivalent to the following: let $\tau = (i_2, \dots, i_d), \tau' = (i_2', \dots, i_d') \in \{\tau \in T: \sum \tau - t < d-1\}$.
    \begin{itemize}
        \item If $I(\tau) = I(\tau')$, then $\tau$ precedes $\tau'$ if and only if $\tau <_{\lex} \tau'$.
        \item If $I(\tau) < I(\tau')$, then let $m$ be the smallest in $[2, d]$ such that $i_m \ne i_m'$. Let $\ell = \min\{m, I(\tau)\}$. In this case, $\tau$ precedes $\tau'$ if and only if $i_\ell \le i_\ell'$. 
    \end{itemize}
\end{lemma}
\begin{proof}
    This follows by construction. 
\end{proof}
Finally, we extend $\prec$ on $T$ to $\{1^1\cup \tau \cup o_k: \tau \in T\}$ and $\{n^1\cup \tau \cup o_k: \tau \in T\}$ respectively and keep the notation $\prec$ for the two sets. 

\subsection{Minimal new faces of the shelling}\label{subsec:minimalnewface}

So far, we have explained the geometric intuition behind the shelling of $\widetilde{B}$ proposed in Table \ref{tab:shellingorder}. In Lemma \ref{lem:paddings}, we characterized all facets of $\partial \widetilde{B}$ that are not contained in $\bigcup_{i_1 \in [n-1]} i_1^1\ i_1+1^1*\widetilde{B(i_1)}$. In Section \ref{subsec:connecting}, we defined an order $\prec$ on $\{\sigma \setminus 2^1 \cup o_k: \sigma \in B_k \setminus (L_k \cup U_k) \cap \st_B 1^1 2^1\}$ and $\{\sigma \setminus n-1^1 \cup o_k: \sigma \in B_k \setminus (L_k \cup U_k) \cap \st_B n-1^1 n^1\}$ for every $k$. 

In this section, we prove that the order of facets of $\widetilde{B}$ in Table \ref{tab:shellingorder} is a shelling of $\widetilde{B}$. Recall the hypothesis of Proposition \ref{prop:mainprop} that $\widetilde{B_i}$ is shellable for every $i_1 \in [n-1]$. Therefore, there is nothing to check for $1^1 2^1 * \widetilde{B(1)}$ at the beginning. 

We focus on the rest of the facets $F \in \widetilde{B}$. By Definition \ref{def:shellability}, it suffices to find for every $F$ an $X(F) \subset F$ such that 
\begin{itemize}
\item $X(F)$ is not contained in any facet preceding $F$ in Table \ref{tab:shellingorder} ($X(F)$ is a new face),
\item for every $x \in X(F)$, $F \setminus x$ is contained in a facet $F'$ preceding $F$ ($X(F)$ is minimal). 
\end{itemize}
Indeed, suppose $G$ is a facet preceding $F$ such that $ \dim (F \cap G) < D - 1$. Then $F \cap G \ne X(F)$ because $F \cap G$ is not a new face. Take any $x \in X(F) \setminus (F \cap G)$. There is a facet $F'$ preceding $F$ such that the $F \setminus x \subseteq F'$. Thus the ridge $F \cap F' = F\setminus x$ contains $F \cap G$. The face $X(F)$ is called the \emph{minimal new face} of $F$ in the shelling. We find the minimal new faces of the facets of $\widetilde{B} \setminus 1^1 2^1 * \widetilde{B(1)}$ in Lemmas \ref{lem:layersminface}-\ref{lem:padding2minface} below.

\begin{lemma}\label{lem:layersminface}
    Let $F$ be a facet of $i_1^1\ i_1+1^1 * \widetilde{B(i_1)}$ such that $i_1 \in [2, n-1]$. Then $X(F) = i_1+1^1 \cup Y(F\setminus i_1^1\ i_1+1^1)$, where $Y(F\setminus i_1^1\ i_1+1^1)$ is the minimal new face of $F\setminus i_1^1\ i_1+1^1$ in the shelling of $\widetilde{B(i_1)}$. 
\end{lemma}
\begin{proof}
    No facet $F'$ of $i_1^1\ i_1+1^1 * \widetilde{B(i_1)}$ can contain $X(F)$, otherwise $Y(F) \subset F'\setminus i_1^1\ i_1+1^1$ contradicts $Y(F)$ being new. Any other facets preceding $F$ miss $i_1+1^1$, so $X(F)$ is a new face.

    If $x \in Y(F\setminus i_1^1\ i_1+1^1)$, then $F \setminus x \subset F'$ for some facet $F'$ such that $F' \setminus i_1^1\ i_1+1^1$ precedes $F \setminus i_1^1\ i_1+1^1$ in the shelling of $\widetilde{B(i_1)}$. It remains to show that $F \setminus i_1+1^1 \subset F'$ for some $F'$ preceding $F$. We gather the case discussions in the table below for clarity. Suppose $\sigma = (i_1, \cdots, i_d)$. 

\setlength{\arrayrulewidth}{0.3mm}
\renewcommand{\arraystretch}{1.6}
\renewcommand\cellset{\renewcommand\arraystretch{1.6}}
\begin{table}[H]
\footnotesize
\begin{center}
\begin{tabular}{ll|ll}
\hline
\multicolumn{1}{l|}{\multirow{3}{*}{$\sigma \in L_k$}} & \multirow{2}{*}{$T_{\sigma, 1}$} & \multicolumn{1}{l|}{$F = \sigma$} & \makecell[l]{Let $\sigma':= (i_1 - 1^1, \dots, i_d^d) \in U_{k-1}$. \\
Take $F' = \sigma'$ for $T_{\sigma', 1}$, \\
$F' = \sigma' \setminus i_1-1^1 \cup o_k$ for $T_{\sigma', 2}$.} \\ \cline{3-4} 
\multicolumn{1}{l|}{}                                       &                     & \multicolumn{1}{l|}{\makecell[l]{$F = \sigma \setminus i_\ell^\ell \cup o_k$,\\$i_\ell \ne n-1$}} & 
\makecell[l]{
Let $\sigma':= (i_1-1^1, \dots, i_\ell+1^\ell, \dots, i_d^d) \in L_k$.\\
Take $F' = \sigma'\setminus i_1-1^1 \cup o_k$ for $T_{\sigma', 1}$, \\
$F' = \sigma' \setminus i_\ell+2^\ell \cup o_k$ for $T_{\sigma', 2}$.
}\\ \cline{2-4} 
\multicolumn{1}{l|}{}                                       & $T_{\sigma, 2}$                  & \multicolumn{2}{l}{Take $F' = \sigma \setminus i_1+1^1 \cup o_k$.}       \\ \hline
\multicolumn{1}{l|}{\multirow{3}{*}{$\sigma \in U_k$}} & $T_{\sigma, 1}$                  & \multicolumn{2}{l}{Take $F' = \sigma \setminus i_1+1^1 \cup o_k$.}       \\ \cline{2-4} 
\multicolumn{1}{l|}{}                                       & \multirow{2}{*}{$T_{\sigma, 2}$} & \multicolumn{1}{l|}{\makecell[l]{$F = \sigma \setminus i_\ell + 1^\ell \cup o_k,$\\ $i_\ell = 1$}} & 
\makecell[l]{
Let $\tau := i_1-1^1\ i_1^1 \cdots i_\ell^\ell \cdots i_d^d\ i_d+1^d \in C_k$.\\
Take $F' = \tau \cup o_k$.
}\\ \cline{3-4} 
\multicolumn{1}{l|}{}                                       &                     & \multicolumn{1}{l|}{$F = \sigma \setminus i_\ell^\ell \cup o_k$} &
\makecell[l]{
If $i_\ell = n - 1$, then let\\
$\tau := i_1-1^1\ i_1^1 \cdots i_\ell+1^\ell \cdots i_d^d\ i_d+1^d \in C_k$. \\
Take $F' = \tau \cup o_k$.\\
If $i_\ell \ne n - 1$, then let\\
$\sigma':= (i_1-1^1, \dots, i_\ell+1^\ell, \dots, i_d^d) \in U_k$.\\
Take $F' = \sigma'\setminus i_\ell+2^\ell \cup o_k$ for $T_{\sigma', 1}$, \\
$F' = \sigma' \setminus i_1-1^1 \cup o_k$ for $T_{\sigma', 2}$.
}\\ \hline
\multirow{2}{*}{$\sigma \in B_k \setminus (L_k \cup U_k)$}                    & \multirow{2}{*}{}   & \multicolumn{1}{l|}{\makecell[l]{$F = \sigma \setminus i_\ell+1^\ell \cup o_k,$\\$ i_\ell = 1$}} & 
\makecell[l]{
Let $\sigma' := (i_1-1^1, \dots, i_d^d)$.\\
If $\sigma \in L_k$, then take $F' = \sigma' \setminus i_1-1^1 \cup o_k$ for $T_{\sigma', 1}$, \\
$F' = \sigma \setminus i_\ell+1^\ell \cup o_k$ for $T_{\sigma', 2}$.\\
Otherwise, take $F' = \sigma' \setminus i_\ell+1^\ell \cup o_k \in C_k *o_k$. 
} \\ \cline{3-4} 
                                       &               & \multicolumn{1}{l|}{\makecell[l]{$F = \sigma \setminus i_\ell^\ell \cup o_k,$\\$i_\ell = n-1$}} & \makecell[l]{
Let $\sigma' := (i_1-1^1, \dots, i_d^d)$.\\
If $\sigma \in L_k$, then take $F' = \sigma' \setminus i_1-1^1 \cup o_k$ for $T_{\sigma', 1}$,\\
$F' = \sigma \setminus i_\ell^\ell \cup o_k$ for $T_{\sigma', 2}$.\\
Otherwise, take $F' = \sigma' \setminus i_\ell^\ell \cup o_k \in C_k *o_k$. 
} \\ \hline
\end{tabular}
\end{center}
\end{table}
\end{proof}

\begin{lemma}\label{lem:padding1minface}
    Let $F$ be a facet in case 1 of Lemma \ref{lem:paddings} such that $F$ comes after $i_1^1\ i_1+1^1 * \widetilde{B(i_1)}$, so that $F$ is either $\sigma \setminus i_1^1 \cup o_k$ with $i_1 \in [n-1]$ or $\sigma \setminus i_1+1^1 \cup o_k$ with $i_1 = 1$ for some $\sigma = (i_1, \dots, i_d) \in L_k \cup U_k \cap \st_B i_1^1\ i_1+1^1$. Then $X(F)$ is the minimal new face of $F$ when we restrict the order in Table \ref{tab:shellingorder} to $T_{\sigma, j}$. 
\end{lemma}
\begin{proof}
     The minimality of $X(F)$ follows from its minimality in the restricted order to $T_{\sigma, j}$. 
     
     We now show that $X(F)$ is a new face. Recall that $T_{\sigma, 1} = F_\sigma * \partial(G_\sigma * o_k)$ and $T_{\sigma, 2} = \partial F_\sigma * (G_\sigma * o_k)$. If $F = \sigma \setminus i_1^1 \cup o_k$, then $T_{\sigma, j} \setminus F$ precedes $F$. So $X(F) = G_\sigma * o_k \setminus i_1^1$ for $T_{\sigma, 1}$ and $X(F) = F_\sigma \setminus i_1^1$ for $T_{\sigma, 2}$. If $F = \sigma \setminus i_1+1^1 \cup o_k$ with $i_1 = 1$, then all facets of $T_{\sigma, j} \setminus F$ except $\sigma \setminus i_1^1 \cup o_k$ precede $F$. So $X(F) = G_\sigma * o_k \setminus i_1^1\ i_1+1^1$ for $T_{\sigma, 1}$ and $X(F) = F_\sigma \setminus i_1^1 i_1+1^1$ for $T_{\sigma, 2}$. 

     Suppose for contradiction that $X(F)$ is contained in a facet $F' \notin T_{\sigma, j}$ preceding $F$. Suppose that $F'$ arises from $\sigma' = (i_1', \dots, i_d') \ne \sigma$. Then by Table \ref{tab:shellingorder}, $i_1' \le i_1$. Using the explicit expressions of $F_\sigma$ from Lemma \ref{lem:missingfaces}, we finish the rest of the proof in the table below. 
\setlength{\arrayrulewidth}{0.3mm}
\renewcommand{\arraystretch}{1.6}
\renewcommand\cellset{\renewcommand\arraystretch{1.6}}
\begin{table}[H]
\footnotesize
\begin{center}
\begin{tabular}{l|l|l|l}
\hline
\multirow{2}{*}{$\sigma \in L_k$} & $T_{\sigma, 1}$ & $G_\sigma = \{i_\ell^\ell: i_\ell \ne n-1\}$ & \makecell[l]{
$F' \supset X(F) = G_\sigma * o_k \setminus i_1^1$\\ implies $\sum \sigma' < \sum \sigma$.\\
But then $\sigma' \notin B_k$, contradicting $o_k \notin F'$. 
} \\ \cline{2-4} 
                     & $T_{\sigma, 2}$ & $F_\sigma = \{i_\ell^\ell: i_\ell = n-1\} \cup \{i_\ell + 1^\ell: \ell \in [d]\}$ & \makecell[l]{
$F' \supset X(F) \in  \{F_\sigma \setminus i_1^1, F_\sigma \setminus i_1^1\ i_1+1^1$\}\\
implies $\sum \sigma' - \sum \sigma \le d-1$.\\
So $\sigma' \in B_k \setminus (L_k \cup U_k)$.\\
$\sigma' \setminus i_1^1 \cup o_k$ comes after $F$,\\
which leaves $F' = \sigma' \setminus i_\ell + 1^\ell \cup o_k$\\
for some $i_\ell + 1 = n-1$. But $i_\ell + 1^\ell \in X(F)$.
} \\ \hline
\multirow{2}{*}{$\sigma \in U_k$}  & $T_{\sigma, 1}$ & $G_\sigma = \{i_\ell + 1^\ell: i_\ell \ne 1\}$ & 
\makecell[l]{
Not applicable\\
since $F = \sigma \setminus i_1 + 1 \cup o_k$\\
comes before $i_1^1\ i_1+1^1 * \widetilde{B(i_1)}$ for all $i_1 \ne 1$. 
} \\ \cline{2-4} 
                     & $T_{\sigma, 2}$ & $F_\sigma = \{i_\ell+1^\ell: i_\ell = 1\} \cup \{i_\ell^\ell: \ell \in [d]\}$ & \makecell[l]{
$F' \supset X(F) \in  \{F_\sigma \setminus i_1^1, F_\sigma \setminus i_1^1\ i_1+1^1$\},\\
implies $\sum \sigma - \sum \sigma' \le d-1$.\\
So $\sigma' \in B_k \setminus (L_k \cup U_k)$.\\
$\sigma' \setminus i_1^1 \cup o_k$ comes after $F$,\\
which leaves $F' = \sigma' \setminus i_\ell^\ell \cup o_k$\\
for some $i_\ell - 1 = 1$. But $i_\ell^\ell \in X(F)$.
} \\ \hline
\end{tabular}
\end{center}
\end{table}
\end{proof}

\begin{lemma}\label{lem:padding2minface}
Let $F$ be a facet in case 1 of Lemma \ref{lem:paddings} such that $F$ comes before $i_1^1\ i_1+1^1 * \widetilde{B(i_1)}$, so that $F = \sigma \setminus i_1+1^1 \cup o_k$ with $i_1 \ne 1$ for some $\sigma = (i_1, \dots, i_d) \in L_k \cup U_k \cap \st_B i_1^1\ i_1+1^1$. Let $Y(F)$ be the minimal face of $F$ that is not in $T_{\sigma, j} \setminus F$. Then $X(F) = (F \setminus Y(F)) \setminus i_1^1$ if $\sigma \in L_k$ and $X(F) = F \setminus Y(F)$ if $\sigma \in U_k$. 
\end{lemma}
\begin{proof}
    The triangulations where $F$ can arise are $T_{\sigma, 2}$ for $\sigma \in L_k$ and $T_{\sigma, 1}$ for $\sigma \in U_k$.

    We first show that $X(F)$ is a new face. Suppose for contradiction that $X(F)$ is contained in a facet $F'$ preceding $F$. Suppose that $F'$ arises from $\sigma' = (i_1', \dots, i_d') \ne \sigma$. 

    If $F \in T_{\sigma, 2}$ for $\sigma \in L_k$, then $X(F) = G_\sigma * o_k \setminus i_1^1 = \{i_\ell^\ell: \ell \in [2, d],\ i_\ell \ne n-1\} \cup o_k$. But this means $\sum \sigma' < \sum \sigma$, so $\sigma' \notin B_k$, contradicting $o_k \in F'$. 

    If $F \in T_{\sigma, 1}$ for $\sigma \in U_k$, then $X(F) = F_\sigma = \{i_\ell^\ell: \ell \in [d]\} \cup \{i_\ell +1^\ell: \ell \in [2,d],\ i_\ell = 1\}$. This means $\sum \sigma - \sum \sigma' \le d$, so $\sigma' \in B_k \setminus (L_k \cup U_k)$. But then $F'$ must be missing either $i_\ell^\ell$ for some $i_\ell = n-1$ or $i_\ell+1^\ell$ for some $i_\ell = 1$. 

    Next, for every $x \in X(F)$, we find a facet $F'$ preceding $F$ containing $F \setminus x$. The argument for the case of $F \in T_{\sigma, 2}$ for $\sigma \in L_k$ is identical to the case of $T_{\sigma, 1}$, $\sigma \in L_k$ in the table of Lemma \ref{lem:layersminface}. The argument for the case of $F \in T_{\sigma, 1}$ for $\sigma \in U_k$ is identical to the case of $T_{\sigma, 2}$, $\sigma \in U_k$ in the table of Lemma \ref{lem:layersminface}. 
\end{proof}

\begin{lemma}\label{lem:connectminface}
    Let $F = a \cup \tau \cup o_k$ be a facet in case 2 of Lemma \ref{lem:paddings} that arises from $\sigma = (i_1, \dots, i_d)$. Recall the definition of $I(\tau)$ from Section \ref{subsec:connecting}. Then
    \[
    X(F) = (a \cap n^1) \cup \bigcup_{\ell \in [2, d]}\{i_\ell^\ell: i_\ell = n - 1 \text{ or } \ell \ge I(\tau)\}  \cup \bigcup_{\ell \in [2, d]}\{i_\ell+1^\ell: i_\ell = 1 \text{ or } \ell \le I(\tau)\}.
    \]
\end{lemma}
\begin{proof}
    Suppose $F'$ is a facet containing $X(F)$ and $F'$ arises from some $\sigma' = (i_1', \dots, i_d')$. Then $X(F) \subset \sigma'$. Because $\sigma \in B_k \setminus (L_k \cup U_k)$, if $\sigma' = \sigma$, then $F'$ misses at least an $i_\ell^\ell$ for some $i_\ell = n-1$ or an $i_\ell + 1^\ell$ for some $i_\ell = 1$. Therefore, $\sigma' \ne \sigma$. Moreover, if $i_1 = 1$ and $i_1' = n-1$, then $F'$ comes after $F$. So assume $i_1 = i_1'$ for the discussion. 

    If $I(\tau) \le 1$, then $X(F) = (a\cap n^1) \cup \{i_\ell: \ell \ge 2\} \cup \{i_\ell+1^\ell: i_\ell = 1, \ell \ge 2\}$. For all $\ell \in [2, d]$, $i_\ell - i_\ell' \in \{0, 1\}$. So $\sum \sigma - \sum \sigma' \le d-1$. This means $\sigma'$ is also in $B_k \setminus (L_k \cup U_k)$, and $F' = a \cup \tau' \cup o_k$. If $I(\tau') \in [2, d]$, then $\tau \prec \tau'$. If $I(\tau') \le 1$, then again $\tau \prec \tau'$ because $\tau' <_{\lex} \tau$. Therefore, $F \prec F'$. 

    For the rest of the proof, suppose $I(\tau) \in [2, d]$. For every $\ell \in [I(\tau)+1, d]$, $i_\ell' - i_\ell \in \{-1, 0\}$. And for every $\ell \in [2, I(\tau) - 1]$, $i_\ell' - i_\ell \in \{0, 1\}$. Thus 
    \[
    \sum \sigma' - \sum \sigma \in [I(\tau) - d, I(\tau) - 2]. 
    \]
    This means $\sigma'\in B_k \setminus (L_k \cup U_k)$, and $F' = a \cup \tau' \cup o_k$ for some $I(\tau') \in [2, d]$. We compare $\tau$ and $\tau'$ with respect to $\prec$.

    First suppose $I(\tau') < I(\tau)$. Note that $I(\tau) - I(\tau')$ is not larger than the number of $\ell \in [2, I(\tau)]$ such that $i_\ell' - i_\ell = 1$. Let $\ell_{\text{min}}$ be the smallest of such $\ell$, then $\ell_{\text{min}} \le I(\tau')$. By the second bullet point of Lemma \ref{lem:connectingalgorithm}, this implies $\tau \prec \tau'$.

    Suppose alternatively that $I(\tau') > I(\tau)$. Let $\ell_{\text{min}} = \min\{\ell \in [2, d]: i_\ell' \ne i_\ell\}$. If $\ell_{\text{min}} < I(\tau)$, then it must be that $i_{\ell_{\text{min}}}' > i_{\ell_{\text{min}}}$. Or $\ell_{\text{min}} > I(\tau)$. In both cases, $\tau \prec \tau'$ follows again from the second bullet point of Lemma \ref{lem:connectingalgorithm}. 

    Finally, suppose $I(\tau') = I(\tau)$. This means there must be some $\ell < I(\tau')$ such that $i_\ell' - i_\ell = 1$, so that $\tau <_{\lex} \tau'$. By the first bullet point of Lemma \ref{lem:connectingalgorithm}, $\tau \prec \tau'$. 

    Therefore, $F'$ always comes after $F$ in the order, proving that $X(F)$ is a new face. To show the minimality of $X(F)$, we find a preceding facet $F'$ containing $F\setminus x$ for each $x \in X(F)$. 

    If $x = a = n^1$, then take $F' = F \setminus a \cup 1^1$. 
    
        If $x = i_\ell^\ell$ for some $\ell \ge 2$, $i_\ell = n-1$, then take $F' = \sigma \setminus i_\ell^\ell \cup o_k \in i_1^1\ i_1+1^1 * \widetilde{B(i_1)}$. 
        
        If $x = i_\ell^\ell$ for some $\ell \ge 2$, $\ell \ge I(\tau)$, and $i_\ell \ne n-1$, then consider $\sigma' = (i_1^1, \dots, i_\ell+1^\ell, \dots, i_d^d)$. If $\sigma' \in U_k$, then take $F' = \sigma' \setminus i_\ell+2^\ell \cup o_k$ for $T_{\sigma', 1}$ and $F' = \sigma' \setminus i_1^1 i_1+1^1 \cup a o_k$ for $T_{\sigma', 2}$. If $\sigma' \in B_k \setminus (L_k \cup U_k)$, then take $F' = \sigma' \setminus i_1^1 i_1+1^1 \cup a o_k$. Since $\sum \sigma' = 1+ \sum \sigma$, $F' = a \cup \tau' \cup o_k$ for some $I(\tau') < I(\tau)$. Since $I(\tau) \le \ell$, $I(\tau') < \ell$. Moreover, $F <_{\lex} F'$. Hence $F' \prec F$ always holds. 

        If $x = i_\ell+1^\ell$ for some $\ell \ge 2$, $i_\ell = 1$, then take $F' = \sigma \setminus i_\ell+1^\ell \cup o_k \in i_1^1\ i_1+1^1 * \widetilde{B(i_1)}$. 

        If $x = i_\ell+1^\ell$ for some $\ell \ge 2$, $\ell \le I(\tau)$, and $i_\ell \ne 1$, then consider $\sigma' = (i_1^1, \dots, i_\ell-1^\ell, \dots, i_d^d)$. If $\sigma' \in L_k$, then take $F' = \sigma' \setminus i_\ell-1^\ell \cup o_k$ for $T_{\sigma', 1}$ and $F' = \sigma' \setminus i_1^1 i_1+1^1 \cup a o_k$ for $T_{\sigma', 2}$. If $\sigma' \in B_k \setminus (L_k \cup U_k)$, then take $F' = \sigma' \setminus i_1^1 i_1+1^1 \cup a o_k$. Since $\sum \sigma' = -1+ \sum \sigma$, $F' = a \cup \tau' \cup o_k$ for some $I(\tau') > I(\tau) \ge \ell \ge 2$. Thus $F' \prec F$. 
\end{proof}

\section{The shellability of the sphere $\widetilde{B} \cup (\partial \widetilde{B} * o)$}
\label{sec:shellablesphere}

\begin{figure}[H]
    \centering
    \includegraphics[width=0.5\linewidth]{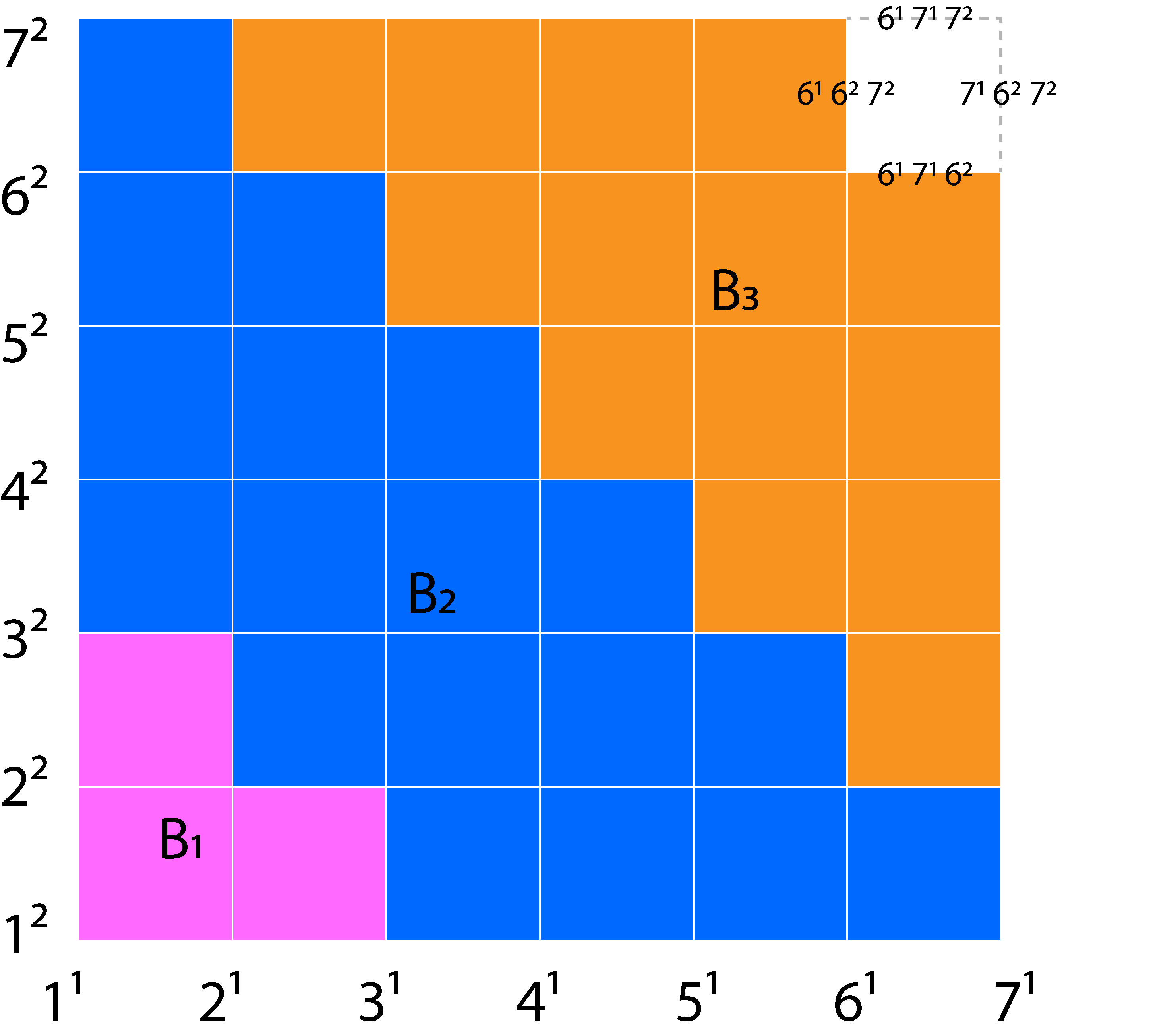}
    \caption{A grid representing $B$ for $d = 2$ and $n = 7$. The facets $6^1 7^1 7^2$ and $7^1 6^2 7^2$ of $\partial B$ are not included in $\partial \widetilde{B}$. The latter boundary contains instead $6^1 6^2 7^2$ and $6^1 7^1 6^2$.}
    \label{fig:boundaryedgecase}
\end{figure}

Observe that $\partial \widetilde{B} = \partial B$ except when $d+2$ divides $d(n-1)$. See Figure \ref{fig:boundaryedgecase} for an example of $d = 2, n = 7$. In this case, the ``uppermost right" cube $(n-1, \dots, n-1) =: \sigma_{\max}$ is not included in any $B_k$. As a result, the ridges $\sigma_{\max}\setminus n-1^m$ for $m \in [d]$ are not facets of $\partial \widetilde{B}$. Instead, $\partial \widetilde{B}$ intersects the uppermost right cube at $\bigcup_{m \in [d]} (\sigma_{\max}\setminus n^m)$. 

Given any face of $B$, we can identify each of its vertices $a^m$ with $a + (m-1)n$. For example, $7^16^27^2$ can be expressed as $7\ 13\ 14$. When we compare two facets of $\partial B$ with respect to lexicographic order below, we mean comparing them using this identification. Moreover, for a facet $\tau$ of $\partial B$, let $\sigma(\tau)$ be the unique facet of $B$ containing $\tau$.

To prove the shellability of $\partial \widetilde{B}$ in both scenarios, we first find a shelling for $\partial B$ that ends with $\{\sigma_{\max}\setminus n-1^m : m \in [d]\}$, then replace it with $\{\sigma_{\max}\setminus n^m : m \in [d]\}_{<_{\lex}}$ and prove that the order remains a shelling. 

Define a total order $\prec_{\partial}$ on the set of facets of $\partial B$ as follows: for any two facets $\tau, \tau'$ of $\partial B$, write $\tau' \prec_\partial \tau$ if and only if
    \begin{itemize}
        \item $\sigma(\tau') = \sigma(\tau)$ and $\tau' <_{\lex} \tau$, or
        \item $\sigma(\tau') <_{\lex} \sigma(\tau)$.
    \end{itemize}
\begin{lemma}\label{lem:boundaryshelling}
    The order $\prec_\partial$ gives a shelling of $\partial B$.
\end{lemma}
\begin{proof}
Let $\tau = i_1^1\ i_1+1^1 \cdots a \cdots i_d^d\ i_d+1^d$ where $a = 1^m$ or $a = n^m$. We claim that 
    \[
    X := (a \cap n^m) \cup \{i_\ell+1^\ell: \ell < m, i_\ell \ne 1\} \cup \{i_\ell+1^\ell: \ell > m\} \cup \{i_\ell^\ell: \ell >m, i_\ell = n-1\} 
    \]
    is a minimal new face of $\tau$ with respect to $\prec_{\partial}$. Any facet $\tau'$ such that $\sigma(\tau') <_{\lex} \sigma(\tau)$ must miss at least one $i_\ell+1^\ell$, so $X \notin \tau'$. On the other hand, if $\sigma(\tau') = \sigma(\tau)$ and $\tau' <_{\lex} \tau$, then $\tau'$ must be $\sigma(\tau) \setminus i_\ell+1^\ell$ where $i_\ell = 1$ or $\sigma(\tau) \setminus i_\ell^\ell$ where $i_\ell = n-1$ for some $\ell > m$. Thus $\tau'$ cannot contain $X$ either. Lastly, for every $x \in X$, we find a $\tau' \prec_\partial \tau$ containing $\tau \setminus x$ in the table below:
\setlength{\arrayrulewidth}{0.3mm}
\renewcommand{\arraystretch}{1.6}
\begin{center}
\begin{tabular}{l|l|l}
\hline
$x$ & $\tau'$ containing $\tau \setminus x$ & Justification\\
\hline
$i_\ell + 1^\ell, \ell \in [d], i_\ell \ne 1$& 
$\tau \setminus i_\ell + 1^\ell \cup i_\ell - 1^\ell$
& $\sigma(\tau') <_{\lex} \sigma(\tau)$
\\ 
$i_\ell + 1^\ell, \ell > m, i_\ell = 1$& 
$\tau \setminus i_\ell + 1^\ell \cup n^m$&
$\sigma(\tau') = \sigma(\tau)$ and $\tau' <_{\lex} \tau$
\\
$i_\ell^\ell: \ell > m, i_\ell = n-1$&
$\tau \setminus i_\ell^\ell \cup n^m$&
$\sigma(\tau') = \sigma(\tau)$ and $\tau' <_{\lex} \tau$
\\
\hline
\end{tabular}
\end{center}
\end{proof}

\begin{proposition}
    The boundary of $\widetilde{B}$ is shellable. 
\end{proposition}
\begin{proof}
    If $d+2$ does not divide $d(n-1)$, then $\partial \widetilde{B} = \partial B$, so it is shellable by Lemma \ref{lem:boundaryshelling}. Suppose $(d+2)$ divides $d(n-1)$. Then
    \[
    \partial \widetilde{B} = \partial B \setminus \bigcup_{m \in [d]} \left(\sigma_{\max}\setminus n-1^m \right) \cup \bigcup_{m \in [d]} \left(\sigma_{\max}\setminus n^m\right). 
    \]
    Observe that $\{\sigma_{\max}\setminus n-1^m : m \in [d]\}$ is terminal with respect to $\prec_\partial$. Remove these facets from the order and attach $\{\sigma_{\max}\setminus n^m: m \in [d]\}_{<_{\lex}}$. 
    
    To see that this is a shelling of $\partial \widetilde{B}$, take $\tau = \sigma_{\max}\setminus n^m$. Then $X:= \{n-1^\ell: \ell \in [d]\} \cup \{n^\ell: \ell > m\}$ is the minimal new face of $\tau$. Indeed, any facet $\tau' \not \subset \sigma_{\max}$ misses at least one $n-1^\ell$. If $\tau' \subset \sigma_{\max}$ and $\{n^\ell: \l > m\} \subset \tau'$, then $\tau <_{\lex} \tau'$. 
    
    Finally, for any $\ell \in [d]$, $\tau \setminus n-1^\ell \cup n-2^m$ is a facet from $\partial B \setminus \bigcup_{m \in [d]} \left(\sigma_{\max}\setminus n-1^m\right)$ containing $\tau \setminus n-1^\ell$. For $\ell > m$, $\tau \setminus n^\ell \cup n^m \subset \sigma_{\max}$ contains $\tau \setminus n^\ell$ and is lexicographically smaller than $\tau$. This proves the minimality of $X$.
\end{proof}

Therefore, $\partial \widetilde{B} * o$ is also shellable. Simply attaching a shelling of $\partial \widetilde{B} * o$ to a shelling of $\widetilde{B}$ gives a shelling of $\widetilde{B} \cup \partial \widetilde{B} * o$. This proves Theorem \ref{thm:mainthm}.

Combining Theorem \ref{thm:mainthm} with the lower bound in (\ref{eqn:Kalai}), we obtain 
\[
s_{\rm{shell}}(D, N) \ge 2^{\Omega(N^{\lceil D/2 \rceil})} \text{ for all } D \ge 3.
\]
On the other hand, we can compute an upper bound for $s_{\rm{shell}}(D, N)$. As mentioned towards the end of Section \ref{sec:intro}, this uses the result of Benedetti and Ziegler \cite{Benedetti} on the number of LC $D$-spheres with $M$ facets along with the Upper Bound Theorem \cite{Stanley} for simplicial spheres: 
\[
s_{\rm{shell}}(D, N) \le \sum_{M = 1}^{O(N^{\lceil D/2 \rceil})} 2^{D^2M} = \frac{2^{D^2}(2^{D^2 O(N^{\lceil D/2 \rceil})}-1)}{2^{D^2}-1} = 2^{O(N^{\lceil D/2 \rceil})}. 
\]
This immediately leads to Corollary \ref{cor:maincorollary}: the number of combinatorially distinct shellable $D$-spheres with $N$ vertices is asympotically given by
\[
s_{\rm{shell}}(D, N) = 2^{\Theta(N^{\lceil D/2 \rceil})} \text{ for all } D \ge 3.
\]

\section*{Acknowledgements}
I am extremely grateful to Isabella Novik for her detailed and invaluable comments on the drafts. I am thankful to the anonymous referee for detecting a mistake in the shelling of the boundary and providing suggestions to improve readability of the paper. This research was partially supported by Graduate Fellowship from NSF grant DMS-2246399. 
\bibliographystyle{abbrv} 
\bibliography{ref}

\end{document}